\documentclass[12pt, twoside, eqno]{article}
\usepackage{latexsym}
\usepackage{amssymb}
\usepackage{amsfonts}
\begin{document}

\noindent {\bf \large Fuzzy sets in $\le$--hypergroupoids}\bigskip

\noindent{\bf Niovi Kehayopulu}\bigskip

\noindent{\bf Abstract.} {\small This paper serves as an example to 
show the way we pass from ordered groupoids (ordered semigroups) to 
ordered hypergroupoids (ordered hypersemigroups), from groupoids 
(semigroups) to hypergroupoids (hypersemigroups). The results on 
semigroups (or on ordered semigroup) can be transferred to 
hypersemigroups (or to ordered hypersemigroups) in the way indicated 
in the present paper.\medskip

}
\section{Introduction and prerequisites}An ordered groupoid (: 
$po$-groupoid) is a nonempty set $S$ endowed with an order ``$\le$" 
and a multiplication ``$\cdot$" such that $a\le b$ implies $ca\le cb$ 
and $ac\le bc$ for every $c\in S$. Let $(S,.,\le)$ be an ordered 
groupoid. A nonempty subset $A$ of $S$ is called a {\it left} (resp. 
{\it right}) ideal of $S$ if (1) $SA\subseteq A$ (resp. $AS\subseteq 
A$) and (2) if $a\in A$ and $S\ni b\le a$, then $b\in A$. A subset 
$A$ of $S$ is called an {\it ideal} of $S$ if it is both a left and a 
right ideal of $S$ [2]. A subgroupoid $F$ of $S$ is called a {\it 
filter} of $S$ if (1) $a,b\in S$ such that $ab\in F$ implies $a\in F$ 
and $b\in F$ and (2) if $a\in F$ and $S\ni b\ge a$, then $b\in F$ 
[1].

Given a set $S$, a fuzzy subset of $S$ (or a fuzzy set in $S$) is, by 
definition, an arbitrary mapping of $S$ into the closed interval 
$[0,1]$ of real numbers (Zadeh). Fuzzy sets in ordered groupoids have 
been first considered in [4], where the following concepts have been 
introduced and studied: A fuzzy subset $f$ of an ordered groupoid 
$(S,.,\le)$ is called a {\it fuzzy left} (resp. {\it fuzzy right}) 
{\it ideal} of $S$ if (1) $x\le y$ implies $f(x)\ge f(y)$ and (2) if 
$f(xy)\ge f(y)$ (resp. $f(xy)\ge f(x)$) for every $x,y\in S$. $f$ is 
called a {\it fuzzy ideal} of $S$ if it is both a fuzzy left ideal 
and a fuzzy right ideal of $S$. A fuzzy subset $f$ of $S$ is called a 
{\it fuzzy filter} of $S$ if (1) $x\le y$ implies $f(x)\le f(y)$ and 
(2) if $f(xy)=\min\{f(x),f(y)\}$ for all $x,y\in S$. A fuzzy subset 
$f$ of a groupoid $S$ is called {\it prime} if 
$f(xy)\le\max\{f(x),f(y)\}$ for all $x,y\in S$. For a groupoid $S$ 
and a fuzzy subset $f$ of $S$, the complement of $f$ is the fuzzy 
subset $f': S\rightarrow [0,1]$ of $S$ defined by $f'(x)=1-f(x)$ for 
all $x\in S$. We have seen in [4] that a nonempty subset $A$ of an 
ordered groupoid $S$ is a left (resp. right) ideal of $S$ if and only 
if its characteristic function $f_A$ is a fuzzy left (resp. right) 
ideal of $S$. A nonempty subset $F$ of an ordered groupoid $S$ is a 
filter of $S$ if and only if the fuzzy subset $f_F$ is a fuzzy filter 
of $S$. A fuzzy subset $f$ of an ordered groupoid $S$ is a fuzzy 
filter of $S$ if and only if the complement $f'$ of $f$ is a fuzzy 
prime ideal of $S$.

In the present paper we examine the results of ordered groupoids 
given in [4] for ordered hypergroupoids. We deal with an 
hypergroupoid $(H,\circ)$ endowed with a relation ``$\le$" (not order 
relation, and  so not compatible with the hyperoperation ``$\circ$" 
in general). Though we could call $\sigma$ that relation and 
$\sigma$--hypergroupoid the hypergroupoid endowed with the relation 
$\sigma$, we will show by ``$\le$" the relation and use the term 
$\le$--hypergroupoid, to emphasize the fact that our results hold for 
ordered hypergroupoids as well. As a consequence,  the results in [4] 
also hold in groupoids endowed with a relation ``$\le$" which is not 
an order in general. Our aim is to show the way we pass from ordered 
groupoids to ordered hypergroupoids.

For a groupoid $(S,.)$ we have one operation corresponding to each 
$(a,b)\in S\times S$ the unique element $ab$ of $S$. For an 
hypergroupoid $H$ we have two ``operations". One of them is the 
``operation" between the elements of $H$ which is called 
hyperoperation as it maps the set $H\times H$ into the set of 
nonempty subsets of $H$ and the other is the operation between the 
nonempty subsets of $H$.
We use the terms left (right) ideal, bi-ideal, quasi-ideal instead of 
left (right) hyperideal, bi-hyperideal, quasi-hyperideal and so on, 
and this is because in this structure there are not two kind of left 
ideals, for example, to distinguish them as left ideal and left 
hyperideal. The left ideal in this structure is that one which 
corresponds to the left ideal of groupoids.\section{Main results}
An {\it hypergroupoid} is a nonempty set $H$ with an hyperoperation 
$$\circ : H\times H \rightarrow {\cal P}^*(H) \mid (a,b) \rightarrow 
a\circ b$$on $H$ and an operation $$* : {\cal P}^*(H)\times {\cal 
P}^*(H) \rightarrow {\cal P}^*(H) \mid (A,B) \rightarrow A*B$$ on 
${\cal P}^*(H)$ (induced by the operation of $H$) such that 
$$A*B=\bigcup\limits_{(a,b) \in\,A\times B} {(a\circ b)}$$ for every 
$A,B\in {\cal P}^*(H)$.\smallskip

\noindent The operation ``$*$" is well defined. Indeed: If $(A,B)\in 
{\cal P}^*(H) \times {\cal P}^*(H)$, then $A*B=\bigcup\limits_{(a,b) 
\in\,A\times B} {(a\circ b)}$. For every $(a,b)\in A\times B$, we 
have $(a,b)\in H\times H$, then $(a\circ b)\in {\cal P}^*(H)$, thus 
we get $A*B\in {\cal P}^*(H)$. If $(A,B),(C,D)\in {\cal P}^*(H)\times 
{\cal P}^*(H)$ such that $(A,B)=(C,D)$, then 
$$A*B=\bigcup\limits_{(a,b) \in\,A\times B} {(a\circ 
b)}=\bigcup\limits_{(a,b) \in\,C\times D} {(a\circ b)}=C*D.$$

As the operation ``$*$" depends on the hyperoperation ``$\circ$", an 
hypergroupoid can be also denoted by $(H,\circ)$ (instead of 
$(H,\circ,*)$).

If $H$ is an hypergroupoid then, for every $x,y\in H$, we 
have$$\{x\}*\{y\}=x\circ y.$$Indeed, 
$\{x\}*\{y\}=\bigcup\limits_{\scriptstyle u \in \{ x\} \hfill\atop
\scriptstyle v \in \{ y\} \hfill} {(u\circ v)}=x\circ y.$\smallskip

The following proposition, though clear, plays an essential role in 
the theory of hypergroupoids.\medskip

\noindent{\bf Proposition 1.} {\it Let $(H,\circ)$ be an 
hypergroupoid, $x\in H$ and $A,B\in {\cal P}^*(H)$. Then we have the 
following:

$(1)$ $x\in A*B$ $\Longleftrightarrow$ $x\in a\circ b$ for some $a\in 
A$, $b\in B$.

$(2)$ If $a\in A$ and $b\in B$, then $a\circ b\subseteq A*B$.} 
\medskip

A nonempty subset $A$ of an hypergroupoid $(H,\circ)$ is called a 
left (resp. right) ideal of $H$ if $H*A\subseteq A$ (resp. 
$A*H\subseteq A$). A subset of $H$ which is both a left ideal and a 
right ideal of $H$ is called an {\it ideal} of $H$. A nonempty subset 
$A$ of $H$ is called a {\it subgroupoid} of $H$ if $A*A\subseteq A$. 
Clearly, every left ideal, right ideal or ideal of $H$ is a 
subgroupoid of $H$. \medskip

\noindent{\bf Definition 2.} By a {\it $\le$--hypergroupoid} we mean 
an hypergroupoid $H$ endowed with a relation denoted by ``$\le$". 
\\We write $b\ge a$ if $a\le b$. \medskip

\noindent{\bf Definition 3.} Let $H$ be a $\le$--hypergroupoid. A 
fuzzy subset $f$ of $H$ is called a {\it fuzzy left ideal} of $H$ if

$(1)$ $x\le y \Rightarrow f(x)\ge f(y)$ and

$(2)$ if $f(x\circ y)\ge f(y)$ for all $x,y\in H$.\\With the property 
(2) we mean the following:

(2) if $x,y\in H$ and $u\in x\circ y$, then $f(u)\ge f(y)$.

A fuzzy subset $f$ of $H$ is called a {\it fuzzy right ideal} of $H$ 
if

$(1)$ $x\le y \Rightarrow f(x)\ge f(y)$ and

$(2)$ if $f(x\circ y)\ge f(x)$ for all $x,y\in H$.\\With the property 
(2) we mean:

(2) if $x,y\in H$ and $u\in x\circ y$, then $f(u)\ge f(x)$.

A fuzzy subset of $H$ is called a {\it fuzzy ideal} of $H$ it is both 
a fuzzy left and a fuzzy right ideal of $H$. As one can easily see, a 
fuzzy subset $f$ of $H$ is a fuzzy ideal of $H$ if and only

(1) $x\le y$ implies $f(x)\ge f(y)$ and

(2) if $f(x\circ y)\ge \max\{f(x),f(y)\} \mbox { for all } x,y\in H$ 
in the sense that$$x,y\in H \mbox { and } u\in x\circ y \mbox { 
implies } f(u)\ge \max\{f(x),f(y)\}.$$

Following Zadeh, any mapping $f : H\rightarrow [0,1]$ of a 
$\le$--hypergroupoid $H$ into the closed interval $[0,1]$ of real 
numbers is called a {\it fuzzy subset} of $H$ (or a {\it fuzzy set} 
in $H$) and $f_A$ (: the characteristic function of $A$) is the 
mapping$$f_A : H \rightarrow \{0,1\} \mid x \rightarrow f_A 
(x)=\left\{ \begin{array}{l}
1\,\,\,\,\,$if$\,\,\,\,x \in A\\
0\,\,\,\,$if$\,\,\,\,x \notin A.
\end{array} \right.$${\bf Definition 4.} Let $H$ be a 
$\le$--hypergroupoid. A nonempty subset $A$ of $H$ is called a {\it 
left} (resp. {\it right}) {\it ideal} of $H$ if

$(1)$ $H*A\subseteq A$ (resp. $A*H\subseteq A)$ and

$(2)$ if $a\in A$ and $H\ni b\le a$, then $b\in A$.\medskip

\noindent{\bf Lemma 5.} {\it Let $(H,\circ)$ be an hypergroupoid. If 
$A$ is a left (resp. right) ideal of H, then for every $h\in H$ and 
every $a\in A$, we have $h\circ a\subseteq A$ (resp. $a\circ 
h\subseteq A$). ``Conversely", if $A$ is a nonempty subset of $H$ 
such that $h\circ a\subseteq A$ (resp. $a\circ h\subseteq A$) for 
every $h\in H$ and every $a\in A$, then the set $A$ is a left (resp. 
right) ideal of H}.\medskip

\noindent{\bf Proposition 6.} {\it Let H be a $\le$--hypergroupoid. 
If $L$ is a left ideal of $H$, then $f_L$ is a fuzzy left ideal of 
$H$. ``Conversely", if $L$ is a nonempty subset of $H$ such that 
$f_L$ is a fuzzy left ideal of H, then L is a left ideal of H}. 
\medskip

\noindent{\bf Proof.} {$\Longrightarrow$.} Let $L$ be a left ideal of 
$H$. By definition, $f_L$ is a fuzzy subset of $H$. Let $x\le y$. If 
$y\not\in L$, then $f_L(y)=0$, so $f_L(x)\ge f_L(y)$. If $y\in L$, 
then $H\ni x\le y\in L$ and, since $L$ is a left ideal of $H$, we 
have $x\in L$. Then $f_L(x)=f_L(y)=1$, so $f_L(x)\ge f_L(y)$.
Let now $x,y\in H$ and $u\in x\circ y$. Then $f_L(u)\ge f_L(y)$. 
Indeed: If $y\in L$ then, by by Proposition 1(2), we have $x\circ 
y\subseteq H*L\subseteq L$, so $u\in L$, then $f_L(y)=f_L(u)=1$. If 
$y\notin L$, then $f_L(y)=0\le f_L(u)$.\\{$\Longleftarrow$.} Let 
$x\in H$ and $y\in L$. Then $x\circ y\subseteq L$. Indeed: Let 
$x\circ y\not\subseteq L$. Then there exists $u\in x\circ y$ such 
that $u\notin L$. Since $u\in x\circ y$, by hypothesis, we have 
$f_L(u)\ge f_L(y)$. Since $u\notin L$, we have $f_L(u)=0$. Since 
$y\in L$, we have $f_L(y)=1$, then $0\ge 1$ which is impossible.
Let now $x\in L$ and $H\ni y\le x$. Then $y\in L$. Indeed: Since 
$f_L$ is a fuzzy left ideal of $H$ and $y\le x$, we have $f_L(y)\ge 
f_L(x)$. Since $x\in L$, $f_L(x)=1$. Then we have $f_L(y)\ge 1$. On 
the other hand, $f_L(y)\le 1$, so we have $f_L(y)=1$, and $y\in L$. 
By Lemma 5, $L$ is a left ideal of $H$. $\hfill\Box$\\In a similar 
way we prove the following:\medskip

\noindent{\bf Proposition 7.} {\it Let H be a $\le$--hypergroupoid. 
If $R$ is a right ideal of $H$, then $f_R$ is a fuzzy right ideal of 
$H$. ``Conversely", if $R$ is a nonempty subset of $H$ such that 
$f_R$ is a fuzzy right ideal of H, then R is a right ideal of H}. 
\medskip

\noindent{\bf Proposition 8.} {\it If H is a $\le$--hypergroupoid, a 
nonempty subset $I$ of $H$ is an ideal of H if and only if $f_I$ is a 
fuzzy ideal of H}.\medskip

Now we introduce the concept of filters and fuzzy filters in 
$\le$--hypergroupoids, and we characterize the filters of 
$\le$--hypergroupoids in terms of fuzzy filters.\medskip

\noindent{\bf Definition 9.} Let $H$ be a $\le$--hypergroupoid. A 
nonempty subset $F$ of $H$ is called a {\it filter} of $H$ if

$(1)$ if $x,y\in F$, then $x\circ y\subseteq F$.

$(2)$ if $x,y\in H$ and $x\circ y\subseteq F$, then $x\in F$ and 
$y\in F$.

$(3)$ if $x,y\in H$, then $x\circ y\subseteq F$ or $(x\circ y)\cap 
F=\emptyset$.

$(4)$ if $x\in F$ and $H\ni y\ge x$, then $y\in F$.\\So a filter of 
$H$ is a subgroupoid of $H$ satisfying the conditions (2)--(4). 
\medskip

\noindent{\bf Remark 10.} Let $H$ be a $\le$--hypergroupoid, $F$ a 
filter of $H$ and $x,y\in H$. The following are equivalent:

$(1)$ $x\circ y\subseteq F$ or $(x\circ y)\cap F=\emptyset$.

$(2)$ if $x\notin F$ or $y\notin F$, then $(x\circ y)\cap 
F=\emptyset$.\\Indeed: $(1)\Longrightarrow (2)$. Let $x\notin F$ or 
$y\notin F$. If $x\circ y\subseteq F$ then, since $F$ is a filter, we 
have $x\in F$ and $y\in F$ which is impossible. Thus we have $x\circ 
y\nsubseteq F$. Then, by (2), $(x\circ y)\cap F=\emptyset$ and (1) is 
satisfied.\\$(2)\Longrightarrow (1)$. Let $x\circ y\nsubseteq F$. If 
$x,y\in F$ then, since $F$ is a filter of $H$, we have $x\circ 
y\subseteq F$ which is impossible. Thus we have $x\notin F$ or 
$y\notin F$. Then, by (2), $(x\circ y)\cap F=\emptyset$, and (1) 
holds true.\medskip

\noindent{\bf Definition 11.} Let $H$ be a $\le$--hypergroupoid. A 
fuzzy subset $f$ of $H$ is called a {\it fuzzy filter} of $H$ if

$(1)$ if $x\le y$ implies $f(x)\le f(y)$ and

$(2)$ if $f(x\circ y)=\min\{f(x),f(y)\}$ for every $x,y\in H$\\in the 
sense that if $x,y\in H$ and $u\in x\circ y$, then 
$f(u)=\min\{f(x),f(y)\}$.\medskip

\noindent{\bf Proposition 12.} {\it Let H be a $\le$--hypergroupoid. 
If $F$ is a filter of $H$, then the fuzzy subset $f_F$ is a fuzzy 
filter of H. ``Conversely", if $F$ is a nonempty subset of $H$ such 
that $f_F$ is a fuzzy filter of H, then F is a filter of H}.\medskip

\noindent{\bf Proof.} $\Longrightarrow$. Let $x\le y$. If $x\notin 
F$, then $f_F(x)=0$, so $f_F(x)\le f_F(y)$. If $x\in F$, then 
$f_F(x)=1$. Since $y\in H$ and $y\ge x\in F$, we have $y\in F$. Then 
$f_F(y)=1$, and $f_F(x)\le f_F(y)$.

Let now $x,y\in H$ and $u\in x\circ y$. Then 
$f_F(u)=\min\{f_F(x),f_F(y)\}$. Indeed:

(a) If $x\circ y\subseteq F$, then $x\in F$ and $y\in F$. Also $u\in 
F$. Then $f_F(x)=f_F(y)=f_F(u)=1$, so $f_F(u)= 
\min\{f_F(x),f_F(y)\}$.

(b) Let $x\circ y\nsubseteq F$. Then $x\notin F$ or $y\notin F$ 
(since $x,y\in F$ implies $x\circ y\subseteq F$, impossible), then 
$f_F(x)=0$ or $f_F(y)=0$, and $\min\{f_F(x),f_F(y)\}=0$. On the other 
hand, since $x\circ y\not\subseteq F$, we have $(x\circ y)\cap 
F=\emptyset$. Since $u\in x\circ y$, we have $u\notin F$. Then 
$f_F(u)=0$, so $f_F(u)=\min\{f_F(x),f_F(y)\}$. \smallskip

\noindent $\Longleftarrow$. Let $x,y\in F$. Then $x\circ y\subseteq 
F$. Indeed: Let $u\in x\circ y$. By hypothesis, we have 
$f_F(u)=\min\{f_F(x),f_F(y)\}$. Since $x,y\in F$, we have 
$f_F(x)=f_F(y)=1$. Then $f_F(u)=1$, and $u\in F$. So $F$ is a 
subgroupoid of $H$. Let $x,y\in F$ such that $x\circ y\subseteq F$. 
Then $x\in F$ and $y\in F$. Indeed: Since $x\circ y\in{\cal P}^*(H)$, 
the set $x\circ y$ is nonempty. Take an element $u\in x\circ y$. 
Since $f_F$ is a fuzzy filter of $H$, we have 
$f_F(u)=\min\{f_F(x),f_F(y)\}.$ Suppose $x\notin F$ or $y\notin F$. 
Then $f_F(x)=0$ or $f_F(y)=0$, $\min\{f_F(x),f_F(y)\}=0$ and 
$f_F(u)=0$. On the other hand, since $u\in x\circ y\subseteq F$, we 
have $f_F(u)=1$. We get a contradiction. Let $x,y\in H$ such that
$x\circ y\not\subseteq F$. Then $(x\circ y)\cap F=\emptyset$. Indeed: 
Let $u\in (x\circ y)\cap F$. Since $u\in x\circ y$, we have 
$f_F(u)=\min\{f_F(x),f_F(y)\}$. If $x\notin F$, then $f_F(x)=0$, then 
$f_F(u)=0$. On the other site, since $u\in F$, we have $f_F(u)=1$ 
which is impossible, so $x\in F$. In a similar way we prove that 
$y\in F$, then $(x\circ y)\subseteq F$ which is impossible. Finally, 
let $x\in F$ and $H\ni y\ge x$. Since $f_F$ is a fuzzy filter of $H$, 
we have $1\ge f_F(y)\ge f_F(x)=1$, then $f_F(y)=1$, so $y\in F$. Thus 
$F$ is a filter of $H$. $\hfill\Box$\medskip

In what follows, for a fuzzy subset $f$ of $S$ we introduce the 
concept of the complement $f'$ of $f$ and prove that $f$ is a fuzzy 
filter of $H$ if and only if $f'$ is a fuzzy prime ideal of 
$H$.\medskip

\noindent{\bf Definition 13.} Let $H$ be an hypergroupoid or 
$\le$--hypergoupoid and $f$ a fuzzy subset of $H$. The fuzzy 
subset$$f' : S \rightarrow [0,1] \mbox { defined by } 
f'(x)=1-f(x)$$is called the {\it complement} of $f$ (in 
$H$).\medskip

We remark the following:

(a) If $x\in H$, then $(f')'(x)=1-f'(x)=f(x)$. Thus we have 
$f'':=(f')'=f$.

(b) $f(x)\le f(y)$ $\Longleftrightarrow$ $f'(x)\ge f'(y)$ $(x,y\in 
H)$.

(c) $f(x)=f(y)$ $\Longleftrightarrow$ $f'(x)=f'(y)$ $(x,y\in 
H)$.\medskip

The Proposition 1 in [3] holds for groupoids and hypergroupoids as 
well and we have the following lemma.\medskip

\noindent{\bf Lemma 14.} {\it Let H be an hypergroupoid, f a fuzzy 
subset of H and $x,y\in H$. Then we 
have$$1-\min\{f(x),f(y)\}=\max\{f'(x),f'(y)\}.$$}{\bf Remark 15.} Let 
$H$ be an hypergroupoid, $f$ a fuzzy subset of $H$ and $x,y\in H$. 
The following are equivalent:

$(1)$ $f(x\circ y)=\min\{f(x),f(y)\}$.

$(2)$ $f'(x\circ y)=\max\{f'(x),f'(y)\}$.\\
Indeed: $(1)\Longrightarrow (2)$. Let $u\in x\circ y$. By (1), we 
have $f(u)=\min\{f(x),f(y)\}$. Then, by Lemma 14, we 
$$f'(u)=1-f(u)=1-\min\{f(x),f(y)\}=\max\{f'(x),f'(y)\},$$and (2) 
holds true.\\
$(2)\Longrightarrow (1)$. Let $u\in x\circ y$. By (2) and Lemma 14, 
we have $$f'(u)=\max\{f'(x),f'(y)\}=1-\min\{f(x),f(y)\}.$$
Then $f(u)=1-f'(u)=\min\{f(x),f(y)\}$, and (1) is satisfied. 
$\hfill\Box$\medskip

\noindent{\bf Definition 16.} Let $H$ be a $\le$--hypergroupoid. A 
fuzzy subset $f$ of $H$ is called {\it fuzzy prime ideal} of $H$ if

(1) $x\le y$ implies $f(x)\ge f(y)$ and

(2) $f(x\circ y)=\max\{f(x),f(y)\} \mbox { for all } x,y\in H$\\that 
is, if $x,y\in H$ and $u\in x\circ y$, then 
$f(u)=\max\{f(x),f(y)\}$.\medskip

Which means that a fuzzy subset $f$ of $H$ is called a fuzzy prime 
ideal of $H$ if it is a prime subset of $H$, that is $f(x\circ y)\le 
\max\{f(x),f(y)\}$ for every $x,y\in H$, and at the same time an 
ideal of $H$.\medskip

\noindent{\bf Proposition 17.} {\it Let H be a $\le$--hypergroupoid 
and f a fuzzy subset of H. Then f is a fuzzy filter of H if and only 
if the complement $f'$ of f is a fuzzy prime ideal of H}.\medskip

\noindent{\bf Proof.} $\Longrightarrow$. Let $x\le y$. Since $f$ is a 
fuzzy filter, we have $f(x)\le f(y)$, then $f'(x)\ge f'(y)$. Let now 
$x,y\in H$ and $u\in x\circ y$. Since $f$ is a fuzzy filter, we have 
$f(u)=\min\{f(x),f(y)\}$. Then, $f'(u)=\max\{f'(x),f'(y)\}$ (cf. also 
the proof of Remark 15), thus $f'$ is a fuzzy prime ideal of $H$. \\
$\Longleftarrow$. Let $x\le y$. Since $f'$ is a fuzzy ideal of $H$, 
we have $f'(x)\ge f'(y)$. Then $f(x)\le f(y)$. Let now $x,y\in H$ and 
$u\in x\circ y$. Since $f'$ is a fuzzy prime ideal of $H$, we have 
$f'(u)=\max\{f'(x),f'(y)\}$, then $f(u)=\min\{f(x),f(y)\}.$ Thus $f$ 
is a fuzzy filter of $H$. $\hfill\Box$
{\small\bigskip

\medskip

\noindent Niovi Kehayopulu, University of Athens, Department of 
Mathematics\\
15784 Panepistimiopolis, Athens, Greece\\email: nkehayop@math.uoa.gr

\end{document}